\documentclass[12pt]{article}
\usepackage{amsfonts}
\usepackage{graphicx}
\usepackage{amsmath}
\usepackage{amssymb}
\usepackage{psfrag}
\title{Moment bounds for IID sequences under sublinear expectations\thanks{First version: Agu. 4, 2009. This is the second version on Apr. 15, 2010.}}

\author{ Feng  Hu\thanks{\small{\it E-mail address}:
hufengqf\symbol{64}163.com (F. Hu).}\\
School of Mathematics
\\ Shandong University\\Jinan 250100, China}

\date{}
\begin{document}
\maketitle

\begin{center}{\bf Abstract}\end{center}
\begin{center}
\parbox{0.8\hsize}{ \parindent=0.5 cm
In this paper, with the notion of independent identically
distributed (IID) random variables under sublinear expectations
introduced by Peng [7-9], we investigate moment bounds for IID
sequences under sublinear expectations. We can obtain a moment
inequality for a sequence of IID random variables under sublinear
expectations. As an application of this inequality, we get the
following result: For any continuous function $\varphi$ satisfying
the growth condition $|\varphi(x)|\leq C(1+|x|^p)$ for some $C>0$,
$p\geq1$ depending on $\varphi$, central limit theorem under
sublinear expectations obtained by Peng [8] still holds.}
\end{center}
\textbf{2000 MR Subject Classification}\ \ 60H10, 60G48

\section{ \textbf {Introduction}}

 In classical probability theory, it is well known that for IID
random variables with $E[X_1]=0$ and $E[|X_1|^r]<\infty$ $(r\geq2)$,
$E[|S_n|^r]=O(n^{\frac{r}{2}})$ holds, and hence
$$\sup\limits_{m\geq0}E[|S_{m+n}-S_m|^r]=O(n^{\frac{r}{2}}).\eqno(1)$$ Bounds of this
kind are potentially useful to obtain limit theorems, especially
strong laws of large numbers, central limit theorems and laws of the
iterated logarithm (see, for example, Serfling [10] and Stout [11],
Chapter 3.7).

Since the paper (Artzner et al. [1]) on coherent risk measures,
people are more and more interested in sublinear expectations (or
more generally, convex expectations, see F\"{o}llmer and Schied [4]
and Frittelli and Rossaza Gianin [5]). By Peng [9], we know that a
sublinear expectation $\hat{E}$ can be represented as the upper
expectation of a subset of linear expectations
$\{E_\theta:\theta\in\Theta\}$, i.e.,
$\hat{E}[\cdot]=\sup\limits_{\theta\in\Theta}E_\theta[\cdot]$. In
most cases, this subset is often treated as an uncertain model of
probabilities $\{P_\theta:\theta\in\Theta\}$ and the notion of
sublinear expectation provides a robust way to measure a risk loss
$X$. In fact, nonlinear expectation theory provides many rich,
flexible and elegant tools.

In this paper, we are interested in
$$\overline{E}[\cdot]=\sup\limits_{P\in{\cal P}}E_P[\cdot],$$ where
${\cal P}$ is a set of probability measures.  The main aim of this
paper is to obtain moment bounds for IID sequences under sublinear
expectations.

This paper is organized as follows: in section 2, we give some
notions and lemmas that are useful in this paper. In section 3, we
give our main results including the proofs.

\section{Preliminaries}

In this section, we introduce some basic notions and lemmas. For a
given set ${\cal P}$ of multiple prior probability measures on
$(\Omega,{\cal F}),$ let ${\cal H}$ be the set of random variables
on $(\Omega,{\cal F}).$

For any $\xi\in{\cal H},$ we define a pair of so-called
maximum-minimum expectations $(\overline{E},\underline{E})$ by
$$\overline{E}[\xi]:=\sup_{P\in{\cal P}}E_P[\xi],\ \ \ \underline{E}[\xi]:=\inf_{P\in{\cal
P}} E_{P}[\xi].$$ Without confusion, here and in the sequel,
$E_{P}[\cdot]$ denotes the classical expectation under probability
measure $P$.

Obviously, $\overline{E}$ is a sublinear expectation in the sense
that

\noindent{\bf Definition 2.1} (see Peng [8, 9]). Let $\Omega$ be a
given set and let ${\cal H}$ be a linear space of real valued
functions defined on $\Omega$. We assume that all constants are in
${\cal H}$ and that $X\in{\cal H}$ implies $|X|\in{\cal H}$. ${\cal
H}$ is considered as the space of our ''random variables''. A
nonlinear expectation $\hat{E}$ on ${\cal H}$ is a functional
$\hat{E}$ : ${\cal H}\mapsto R$ satisfying the following
properties: for all $X$, $Y\in{\cal H}$, we have\\
(a) Monotonicity:\ \ \ \ If $X\geq Y$ then
$\hat{E}[X]\geq\hat{E}[Y]$.\\
(b) Constant preserving: $\hat{E}[c]=c$.\\
The triple $(\Omega,{\cal H},\hat{E})$ is called a nonlinear
expectation space (compare with a probability space $(\Omega,{\cal
F},P)$). We are mainly concerned with sublinear expectation where
the expectation $\hat{E}$ satisfies
also\\
(c) Sub-additivity:\ \ \ \ $\hat{E}[X]-\hat{E}[Y]\leq
\hat{E}[X-Y]$.\\
(d)Positive homogeneity: $\hat{E}[\lambda X]=\lambda\hat{E}[X]$,
$\forall\lambda\geq0$.\\
If only (c) and (d) are satisfied, $\hat{E}$ is called a sublinear
functional.

The following representation theorem for sublinear expectations is
very
useful (see Peng [9] for the proof).\\
{\bf Lemma 2.1.} Let $\hat{E}$ be a sublinear functional defined on
$(\Omega,{\cal H})$, i.e., (c) and (d) hold for $\hat{E}$. Then
there exists a family $\{E_\theta:\theta\in\Theta\}$ of linear
functionals on $(\Omega,{\cal H})$ such that
$$\hat{E}[X]=\max\limits_{\theta\in\Theta}E_\theta[X].\eqno(2)$$ If (a) and
(b) also hold, then $E_\theta$ are linear expectations for
$\theta\in\Theta$. If we make furthermore the following assumption:
(H) For each sequence $\{X_n\}_{n=1}^\infty\subset{\cal H}$ such
that $X_n(\omega)\downarrow0$ for $\omega$, we have
$\hat{E}[X_n]\downarrow0$. Then for each $\theta\in\Theta$, there
exists a unique ($\sigma$-additive) probability measure $P_\theta$
defined on $(\Omega,\sigma({\cal H}))$ such that
$$E_\theta[X]=\int_\Omega X(\omega){\rm d}P_\theta(\omega),\ \ X\in{\cal
H}.\eqno(3)$$ {\bf Remark 2.1.} Lemma 2.1 shows that in most cases,
a sublinear expectation indeed is a supremum expectation. That is,
if $\hat{E}$ is  a sublinear expectation on ${\cal H}$ satisfying
(H), then there exists a set (say $\hat{\cal P}$) of probability
measures such that
$$\hat{E}[\xi]=\sup_{P\in\hat{\cal P}}E_P[\xi],\ \ \
-\hat{E}[-\xi]=\inf_{P\in\hat{\cal P}} E_{P}[\xi].$$ Therefore,
without confusion, we sometimes call supremum expectations as
sublinear expectations.

Moreover, a supremum expectation $\overline{E}$ can generate a pair
$(V,v)$ of capacities denoted by
$$V(A):=\overline{E}[I_A],\ \ \ v(A):=-\overline{E}[-I_A],\ \ \forall A\in{\cal F}.$$
It is easy to check that the pair of capacities satisfies
$$V(A)+v(A^c)=1,\ \ \ \forall A\in{\cal F}$$
where $A^c$ is the complement set of $A$.

The following is the notion of IID random variables under sublinear
expectations  introduced by Peng [7-9].

\noindent{\bf Definition 2.2} ({\bf IID under sublinear
expectations}). {\bf Independence:} Suppose that
$Y_1,Y_2,\cdots,Y_n$ is a sequence of random variables such that
$Y_i \in{\cal H}$. Random variable $Y_n$ is said to be independent
of $X:=(Y_1,\cdots,Y_{n-1})$ under $\overline{E}$, if for each
measurable function $\varphi$ on $R^n$ with $\varphi(X,Y_n)\in{\cal
H}$ and $\varphi(x,Y_n)\in{\cal H}$ for each $x\in{R}^{n-1},$ we
have
$$\overline{E}[\varphi(X,Y_n)]=\overline{E}[\overline{\varphi}(X)],$$
where $\overline{\varphi}(x):=\overline{E}[\varphi(x,Y_n)]$ and
$\overline{\varphi}(X)\in{\cal H}$.

{\bf Identical distribution:} Random variables $X$ and $Y$ are said
to be identically distributed, denoted by $X\sim Y$, if for each
measurable function $\varphi$ such that $\varphi(X), \;
\varphi(Y)\in{\cal H}$,
$$\overline{E}[\varphi(X)]=\overline{E}[\varphi(Y)].$$

{\bf IID random variables:} A sequence of random variables
$\{X_i\}_{i=1}^\infty$ is said to be IID, if $X_i\sim X_1$ and
$X_{i+1}$ is independent of $Y:=(X_1,\cdots,X_i)$ for each $i\ge 1.$

\noindent{\bf Definition 2.3} ({\bf Pairwise independence}, see
Marinacci [6]). Random variable $X$ is said to be pairwise
independent of $Y$ under capacity $\hat{V},$ if for all subsets $D$
and $G\in{\cal B}(R),$
$$\hat{V}(X\in D, Y\in G)=\hat{V}(X\in D)\hat{V}(Y\in G).$$

The following lemma shows the relation between Peng's independence
and pairwise independence.

\noindent{\bf Lemma 2.2.} Suppose that $X,Y \in{\cal H}$ are two
random variables. $\overline{E}$ is a sublinear expectation and
$(V,v)$ is the pair of capacities generated by $\overline{E}$. If
random variable $X$ is independent of $Y$ under $\overline{E}$, then
$X$ also is pairwise independent of $Y$ under capacities $V$ and
$v$.

\noindent{\it Proof.} If we choose $\varphi(x,y)=I_D(x)I_G(y)$ for
$\overline{E}$, by the definition of Peng's independence, it is easy
to obtain
$$V(X\in D, Y\in G)=V(X\in D)V(Y\in G).$$

Similarly,  if we choose $\varphi(x,y)=-I_D(x)I_G(y)$ for
$\overline{E}$, it is easy to obtain
$$v(X\in D, Y\in G)=v(X\in D)v(Y\in G).$$
The proof is complete.

Let $C_b(R^n)$ denote the space of bounded and continuous functions,
let $C_{l,Lip}(R^n)$ denote the space of functions $\varphi$
satisfying
$$|\varphi(x)-\varphi(y)|\leq C(1+|x|^m+|y|^m)|x-y|\ \ \ \forall
x,y\in R^n,$$ for some $C>0$, $m\in N$ depending on $\varphi$ and
let $C_{b,Lip}(R^n)$ denote the space of bounded functions $\varphi$
satisfying
$$|\varphi(x)-\varphi(y)|\leq C|x-y|\ \ \
\forall x,y\in R^n,$$ for some $C>0$ depending on $\varphi$.

From now on, we consider the following sublinear expectation space
$(\Omega,{\cal H},\overline{E})$: if $X_1,\cdots,X_n\in{\cal H}$,
then $\varphi(X_1,\cdots,X_n)\in{\cal H}$ for each $\varphi\in
C_{l,Lip}(R^n)$.

\noindent{\bf Definition 2.4} ({\bf $G$-normal distribution}, see
Definition 10 in Peng [7]). A random variable $\xi\in{\cal H}$ under
sublinear expectation $\widetilde{E}$ with
$\overline{\sigma}^2=\widetilde{E}[\xi^2]$,
$\underline{\sigma}^2=-\widetilde{E}[-\xi^2]$ is called $G$-normal
distribution, denoted by ${\cal
N}(0;[{\underline\sigma}^2,{\overline\sigma}^2])$, if for any
function $\varphi\in C_{l,Lip}(R)$, write
$u(t,x):=\widetilde{E}[\varphi(x+\sqrt{t}\xi)],$
$(t,x)\in[0,\infty)\times R$, then $u$ is the unique viscosity
solution of PDE:
$$\partial_t u -G(\partial^2_{xx} u)=0,\ \ u(0,x)=\varphi(x),$$
where $G(x):=\frac{1}{2}(\overline \sigma^2 x^+-\underline \sigma^2
x^-)$ and $x^+:=\max\{x,0\}$, $x^-:=(-x)^+$.

With the notion of IID under sublinear expectations, Peng shows
central limit theorem under sublinear expectations (see Theorem 5.1
in Peng [8]).

\noindent{\bf Lemma 2.3} ({\bf Central limit theorem under sublinear
expectations}). Let $\{X_i\}_{i=1}^\infty$ be a sequence of IID
random variables. We further assume that
$\overline{E}[X_1]=\overline{E}[-X_1]=0.$ Then the sequence
$\{\overline{S}_n\}_{n=1}^\infty$
 defined by $\overline{S}_n:=\frac{1}{\sqrt{n}}\sum\limits_{i=1}^nX_i$ converges in law to $\xi$, i.e.,
$$\lim\limits_{n\rightarrow\infty}\overline{E}[\varphi(\overline{S}_n)]=\widetilde{E}[\varphi(\xi)],$$
for any continuous function $\varphi$ satisfying linear growth
condition (i.e., $|\varphi(x)|\leq C(1+|x|)$ for some $C>0$
depending on $\varphi$), where $\xi$ is a $G$-normal distribution.

\section{ Main results and proofs}

\noindent{\bf Theorem 3.1.} Let a random sequence
$\{X_n\}_{n=1}^\infty$ be IID under $\overline{E}$. Denote
$S_n:=\sum\limits_{i=1}^nX_i$. Assume that
$\overline{E}[X_1]=\overline{E}[-X_1]=0$. Then for each $r>2$, there
exists a positive constant $K_r$ not depending on $n$ such that for
all $n\in N$,
$$\sup\limits_{m\geq0}\overline{E}[|S_{m+n}-S_m|^r]\leq
K_rn^\frac{r}{2}.$$

\noindent{\it Proof.}  Let $r=\theta+\gamma$, where $\theta\in N,
\theta\geq2$ and $\gamma\in(0,1]$. For simplicity, write
$$S_{m,n}:=S_{m+n}-S_m,$$
$$a_n:=\sup\limits_{m\geq0}\overline{E}[|S_{m,n}|^r].$$

Firstly, we shall show that there exists a positive constant $C_r$
not depending on $n$ such that for all $n\in N$,
$$\overline{E}[|S_{m,2n}|^r]\leq2a_n+C_ra_n^{1-\gamma}n^{\frac{\gamma
r}{2}}.\eqno(4)$$ In order to prove (4), we show the following
inequalities for all $n\in N$:
$$\overline{E}[|S_{m,2n}|^r]\leq2a_n+2^{\theta+1}(\overline{E}[|S_{m,n}|^\gamma|S_{m+n,n}|^\theta]
+\overline{E}[|S_{m,n}|^\theta|S_{m+n,n}|^\gamma]),\eqno(5)$$
$$\overline{E}[|S_{m,n}|^\gamma|S_{m+n,n}|^\theta]\leq a_n^{1-\gamma}(\overline{E}[|S_{m,n}|
|S_{m+n,n}|^{\theta-1+\gamma}])^\gamma,\eqno(6)$$
$$\overline{E}[|S_{m,n}|^\theta|S_{m+n,n}|^\gamma]\leq a_n^{1-\gamma}(\overline{E}[|S_{m,n}|^{\theta-1+\gamma}
|S_{m+n,n}|])^\gamma,\eqno(6^{'})$$ $$\overline{E}[|S_{m,n}|
|S_{m+n,n}|^{\theta-1+\gamma}]\leq D_rn^\frac{r}{2},\eqno(7)$$
$$\overline{E}[|S_{m,n}|^{\theta-1+\gamma} |S_{m+n,n}|]\leq D_rn^\frac{r}{2},\eqno(7^{'})$$
where $ D_r$ is a positive constant not depending on $n$.

To prove (5). Elementary estimates yield the following inequality
(*):
$$\begin{array}{lcl}&&|S_{m,2n}|^r=|S_{m,n}+S_{m+n,n}|^{\theta+\gamma}\leq(|S_{m,n}|+|S_{m+n,n}|)^\theta(|S_{m,n}|+|S_{m+n,n}|)^\gamma\\
&\leq&\sum_{i=0}^\theta C_\theta^i|S_{m,n}|^{\theta-i}|S_{m+n,n}|^i(|S_{m,n}|^\gamma+|S_{m+n,n}|^\gamma)\\
&\leq&|S_{m,n}|^{\theta+\gamma}+|S_{m+n,n}|^{\theta+\gamma}+2\sum_{i=0}^\theta
C_\theta^i(|S_{m,n}|^{\gamma}|S_{m+n,n}|^\theta+|S_{m,n}|^{\theta}|S_{m+n,n}|^\gamma)\\
&\leq&|S_{m,n}|^{\theta+\gamma}+|S_{m+n,n}|^{\theta+\gamma}+2^{\theta+1}(|S_{m,n}|
^{\gamma}|S_{m+n,n}|^\theta+|S_{m,n}|^{\theta}|S_{m+n,n}|^\gamma).\end{array}$$
Since $\{X_n\}_{n=1}^\infty$ is a IID random sequence, by the
definition of IID under sublinear expectations,
$$a_n=\sup\limits_{m\geq0}\overline{E}[|S_{m,n}|^r]=\sup\limits_{m\geq0}\overline{E}[|S_{m+n,n}|^r].$$  Taking
$\overline{E}[\cdot]$ on both sides of (*), we have
$$\overline{E}[|S_{m,2n}|^r]\leq2a_n+2^{\theta+1}(\overline{E}[|S_{m,n}|^\gamma|S_{m+n,n}|^\theta]
+\overline{E}[|S_{m,n}|^\theta|S_{m+n,n}|^\gamma]).$$ Hence, (5)
holds.

Since the proof of ($6^{'}$) is very similar to that of (6), we only
prove (6). Without loss of generality, we assume $\gamma\in(0,1)$.
By H\"{o}lder's inequality,
$$\begin{array}{lcl}\overline{E}[|S_{m,n}|^\gamma|S_{m+n,n}|^\theta]&\leq&(\overline{E}[|S_{m,n}||S_{m+n,n}|^{\theta-1+\gamma}])^\gamma(\overline{E}
[|S_{m+n,n}|^{\frac{\theta-\gamma(\theta-1+\gamma)}{1-\gamma}}])^{1-\gamma}\\
&\leq&
a_n^{1-\gamma}(\overline{E}[|S_{m,n}||S_{m+n,n}|^{\theta-1+\gamma}])^\gamma.\end{array}$$
This proves (6).

To prove (7). By the definition of IID under sublinear expectations
and Schwarz's inequality, we have
$$\overline{E}[|S_{m,n}||S_{m+n,n}|^{\theta-1+\gamma}]=\overline{E}[|S_{m,n}|]
\overline{E}[|S_{m+n,n}|^{\theta-1+\gamma}]\leq(\overline{E}[|S_{m,n}|^2])^\frac{1}{2}\overline{E}[|S_{m+n,n}|^{\theta-1+\gamma}].\eqno(8)$$

Next we prove $$\overline{E}[S_{m,n}^2]\leq n\overline{E}[X_1^2],\ \
\forall m\geq0.$$ Indeed, using the definition of IID under
sublinear expectations again, we have
$$\begin{array}{lcl}
&&\overline{E}[S_{m,n}^2]=\overline{E}[(S_{m,n-1}+X_{m+n})^2]=\overline{E}[S_{m,n-1}^2+2S_{m,n-1}X_{m+n}+X_{m+n}^2]\\
&\leq&\overline{E}[S_{m,n-1}^2]+\overline{E}[X_{m+n}^2]\leq\cdots=n\overline{E}[X_1^2].\end{array}$$
So $$\overline{E}[S_{m,n}^2]\leq n\overline{E}[X_1^2]\eqno(9)$$ and
$$\overline{E}[S_{m+n,n}^2]\leq n\overline{E}[X_1^2]\eqno(10)$$ hold. On the other hand, by H\"{o}lder's inequality,
$$\overline{E}[|S_{m+n,n}|^{1+\gamma}]\leq(\overline{E}[S_{m+n,n}^2])^
\frac{1+\gamma}{2}\leq
n^\frac{1+\gamma}{2}(\overline{E}[X_1^2])^\frac{1+\gamma}{2}.\eqno(11)$$
If $\theta=2$, (7) follows from (8), (9), (10) and (11). If
$\theta>2$, we inductively assume
$$\overline{E}[|S_{m+n,n}|^{\theta-1+\gamma}]\leq
M_rn^{\frac{\theta-1+\gamma}{2}},\eqno(12)$$ where $ M_r$ is a
positive constant not depending on $n$. Then (8), (9) and (12) yield
(7). In a similar manner, we can prove that ($7^{'}$) holds.

From (5)-($7^{'}$), it is easy to check that (4) holds. From (4), we
can obtain that for all $n\in N$,
$$a_{2n}\leq2a_n+C_ra_n^{1-\gamma}n^{\frac{\gamma
r}{2}}.$$ By induction, there exists a positive constant $C_r^{'}$
not depending on $n$ such that $a_n\leq C_r^{'}n^\frac{r}{2}$ for
all $n\in\{2^k:k\in N\bigcup\{0\}\}$.

If $n$ is any positive integer, it can be written in the form
$$n=2^k+v_12^{k-1}+\cdots+v_k\leq2^k+2^{k-1}+\cdots+1$$ where
$2^k\leq n<2^{k+1}$ and each $v_j$ is either $0$ or $1$. Then
$S_{m,n}$ can be written as the sum of $k+1$ groups of sums
containing $2^k,v_12^{k-1},\cdots$ terms and using Minkowski's
inequality,
$$\begin{array}{lcl} &&a_n\leq\sup\limits_{m\geq0}[(\overline{E}[|S_{m+v_k+\cdots+v_12^{k-1},2^k}|^r])^\frac{1}{r}+
\cdots+(\overline{E}[|S_{m,v_k}|^r])^\frac{1}{r}]^r\\
&\leq& C_r^{'}[2^\frac{k}{2}+\cdots+1]^r
=C_r^{'}[\frac{2^{\frac{k+1}{2}}-1}{2^{\frac{1}{2}}-1}]^r\leq
K_rn^{\frac{r}{2}}.\end{array}$$
The proof is complete.\\
{\bf Remark 3.1.} {\rm (i)} From the proof of Theorem 3.1, we can
check that the assumption of IID under $\overline{E}$ can be
replaced by the weaker assumption that $\{X_n\}_{n=1}^\infty$ is a
IID random sequence under $\overline{E}$ with respect to the
following functions
$$\varphi_1(x)=x;\ \ \ \varphi_2(x)=-x;$$
$$\varphi_3(x_1,\cdots,x_n)=|x_1+\cdots+x_n|^r,\ \ \ n=1,2,\cdots,\ \ \ r\geq2;$$
$$\varphi_4(x_1,\cdots,x_m,x_{m+1},\cdots,x_{m+n})
=|x_1+\cdots+x_m||x_{m+1}+\cdots+x_{m+n}|^p,\\
m,n=1,2,\cdots,\ \ \ p>1;$$ and
$$\varphi_5(x_1,\cdots,x_m,x_{m+1},\cdots,x_{m+n})=|x_1+\cdots+x_m|^p|x_{m+1}+\cdots+x_{m+n}|,\\
m,n=1,2,\cdots,\ \ \ p>1.$$ {\rm (ii)} A close inspection of the
proof of Theorem 3.1 reveals that the definition of IID under
sublinear expectations plays an important role in the proof. The
proof of Theorem 3.1 is very similar to the classical arguments,
e.g., in Theorem 1 of Birkel [2].

Applying Theorem 3.1, we can obtain the following result:

\noindent{\bf Theorem 3.2.} Let $\{X_i\}_{i=1}^\infty$ be a sequence
of IID random variables. We further assume that
$\overline{E}[X_1]=\overline{E}[-X_1]=0.$ Then the sequence
$\{\overline{S}_n\}_{n=1}^\infty$
 defined by $\overline{S}_n:=\frac{1}{\sqrt{n}}\sum\limits_{i=1}^nX_i$ converges in law to $\xi$, i.e.,
$$\lim\limits_{n\rightarrow\infty}\overline{E}[\varphi(\overline{S}_n)]=\widetilde{E}[\varphi(\xi)],\eqno(13)$$
for any continuous function $\varphi$ satisfying the growth
condition $|\varphi(x)|\leq C(1+|x|^p)$ for some $C>0$, $p\geq1$
depending on $\varphi$, where $\xi$ is a $G$-normal distribution.

\noindent{\it Proof.} Indeed, we only need to prove that (13) holds
for the $p>1$ cases. Let $\varphi$ be an arbitrary continuous
function with growth condition $|\varphi(x)|\leq C(1+|x|^p)$
($p>1$). For each $N>0$, we can find two continuous functions
$\varphi_1$, $\varphi_2$ such that $\varphi=\varphi_1+\varphi_2$,
where $\varphi_1$ has a compact support and $\varphi_2(x)=0$ for
$|x|\leq N$, and $|\varphi_2(x)|\leq|\varphi(x)|$ for all $x$. It is
clear that $\varphi_1\in C_b(R)$ and
$$|\varphi_2(x)|\leq\frac{2C(1+|x|^{p+1})}{N},\ \ \ \hbox{for}\ \ x\in R.$$ Thus
$$\begin{array}{lcl}
|\overline{E}[\varphi(\overline{S}_n)]-\widetilde{E}[\varphi(\xi)]|&=&
|\overline{E}[\varphi_1(\overline{S}_n)+\varphi_2(\overline{S}_n)]-\widetilde{E}[\varphi_1(\xi)+\varphi_2(\xi)]|\\
&\leq&|\overline{E}[\varphi_1(\overline{S}_n)]-\widetilde{E}[\varphi_1(\xi)]|+
|\overline{E}[\varphi_2(\overline{S}_n)]-\widetilde{E}[\varphi_2(\xi)]|\\
&\leq&|\overline{E}[\varphi_1(\overline{S}_n)]-\widetilde{E}[\varphi_1(\xi)]|+
\frac{2C}{N}(2+\overline{E}[|\overline{S}_n|^{p+1}]+\widetilde{E}[|\xi|^{p+1}]).\end{array}$$
Applying Theorem 3.1, we have
$\sup\limits_n\overline{E}[|\overline{S}_n|^{p+1}]<\infty$. So the
above inequality can be rewritten as
$$|\overline{E}[\varphi(\overline{S}_n)]-\widetilde{E}[\varphi(\xi)]|
\leq|\overline{E}[\varphi_1(\overline{S}_n)]-\widetilde{E}[\varphi_1(\xi)]|+\frac{\overline{C}}{N},$$
where
$\overline{C}=2C(2+\sup\limits_n\overline{E}[|\overline{S}_n|^{p+1}]+\widetilde{E}[|\xi|^{p+1}])$.
From Lemma 2.3, we know that (13) holds for any $\varphi\in C_b(R)$
with a compact support. Thus, we have
$\limsup\limits_{n\rightarrow\infty}|\overline{E}[\varphi(\overline{S}_n)]
-\widetilde{E}[\varphi(\xi)]|\leq\frac{\overline{C}}{N}$. Since $N$
can be arbitrarily large, $\overline{E}[\varphi(\overline{S}_n)]$
must converge to $\widetilde{E}[\varphi(\xi)]$. The proof of Theorem
3.2 is complete.

\end{document}